\input amssym.tex
%%%%%%%%%%%%%%%%%%%%%%%
\font\twbf=cmbx12
\font\ninerm=cmr9

\font\sc=cmcsc10

\def\rond{{\scriptstyle\circ}}

\def\cf{{\it cf.\/}\ }    
\def\ie{{\it i.e.\/}\ }
 
\def\up#1{\raise 1ex\hbox{\sevenrm#1}}

\def\cqfd{\hfill\hbox{\vrule\vbox to 4pt{\hrule width
4pt\vfill \hrule}\vrule}}

\def \bg {\bigskip \goodbreak}
\def \sn {\nobreak \smallskip}

\def\ref#1&#2&#3&#4&#5\par{\par{\leftskip = 5em {\noindent
\kern-5em\vbox{\hrule height0pt depth0pt width
5em\hbox{\bf[\kern2pt#1\unskip\kern2pt]\enspace}}\kern0pt}
{\sc\ignorespaces#2\unskip},\
{\rm\ignorespaces#3\unskip}\
{\sl\ignorespaces#4\unskip\/}\
{\rm\ignorespaces#5\unskip}\par}}

\def \eps {\varepsilon}

\def \comp{ \Bbb{C} }
\def \nat{ \Bbb{N}}
 
\def\ent{ \Bbb{Z}}

\def \abs#1{\left\vert#1\right\vert }

\def\N#1{\muskip0=-2mu{\left|\mkern\muskip0\left|
#1\right|\mkern\muskip0\right|}}

\def\pe{\Bbb{P}}
\def\ee{\Bbb{E}}

\def\tore{\Bbb{T}}

\def\wcheck#1{\smash{
        \mathop{#1}\limits^{\scriptscriptstyle{\vee}}}}

\def\Otimes_#1^#2{\matrix{{}_{#2}\cr
      \bigotimes \cr {}^{#1}}}

\def\wco#1#2{\lower0.5ex\hbox{$\scriptstyle{{#1}\wcheck\otimes
{#2}}$}}

\def\refl#1&#2&#3&#4&#5\par{\par{\leftskip = 7em {\noindent
\kern-7em\vbox{\hrule height0pt depth0pt width
5em\hbox{\bf[\kern2pt#1\unskip\kern2pt]\enspace}}\kern0pt}
{\sc\ignorespaces#2\unskip},\
{\rm\ignorespaces#3\unskip}\
{\sl\ignorespaces#4\unskip\/}\
{\rm\ignorespaces#5\unskip}\par}}

\def\ref#1&#2&#3&#4&#5\par{
\item{[{\bf\ignorespaces#1\unskip}]}
{\sc\ignorespaces#2\unskip},\
{\rm\ignorespaces#3\unskip}\
{\sl\ignorespaces#4\unskip\/}\
{\rm\ignorespaces#5\unskip}\par}

\def\liminf_#1{\mathop{\underline{\hbox{\rm lim}}}\limits_{#1}}
\def\limsup_#1{\mathop{\overline{\hbox{\rm lim}}}\limits_{#1}}

\def \ds{\displaystyle}

%%%%%%%%%%%%%%%%%%%%%%%
\centerline{\twbf H$^{\sevenbf 1}$-Projective
Banach Spaces}
\bg
\centerline{{\sc Omran Kouba}}\par
\centerline{Department of Mathematics}\par
\centerline{{\sl Higher Institute for Applied Sciences and Technology}}\par
\centerline{P.O. Box 31983, Damascus, Syria.}\par
\centerline{{\it E-mail} : omran\_kouba@hiast.edu.sy}\par
\bg
{\bf {Abstract:}}
We study the $H^1$-projective
spaces . 
We prove that they have the  Analytic Radon-Nikodym
Property ,
and that they are cotype 2 spaces which satisfy
Grothendieck's Theorem . We show also that the 
ultraproduct of $H^1$-projective spaces is
$H^1$-projective. Other results are also discussed.
\bg
\noindent\bf{1. Introduction, preliminary notation and definitions}\par
\bg
\rm
Let ${\cal D} = \{ z \in \comp : \abs{z} < 1\}$ be the open unit disc
 and let
$dm$ 
be  the normalized Haar measure on its
bounday $\tore$.
 Let $X$ be any complex Banach space and $(\Omega,\cal
A,\pe)$ be any measure  space. 
For $p \in [ 1,+\infty] ~ $ we  denote by $ L^p(\Omega;X) $ ( or 
$ L^p(\Omega,{\cal A},\pe;X) $ if necessary) 
the space of all strongly $\cal A$-measurable functions
$f:\Omega \longrightarrow X $ such that
$\N{f}_{L^p(X)} < +\infty $, with
$$ \N{f}_{L^p(X)} = \left (\int_\Omega \N{f(\omega)}_X^p\,d\pe
(\omega)\right )^{ 1/ p} $$
and with the habitual changes in the case $p =+\infty$.
\par
If $I$ is a set then, $\ell^p(I)$ denotes the Banach
space $L^p(I,\nu)$,where $\nu$ is the counting
measure. And if 
 $ I $ is the set $\{ 1,...,n \} $,  the
corresponding  space  is  denoted  $ \ell_n^p $ .
\par
\bg
We will denote by $ H^p(X)$the space of all analytic functions
$f:{\cal D}\longrightarrow X $ such that $\N{f}_{H^p(X)} < +\infty $,
\
with
$$\N{f}_{H^p(X)} =\hbox{
\
Sup
\
}\{\N{f_r}_{L^p(\tore;X)}:0<r<1\}$$
and $f_r(t) = f(rt) $, for all $t\in\tore$ 
\
and all $r\in[0,1[$.
\par
\bg
We will also denote by $\widetilde H^p(X)$ the closure in
$L^p(\tore;X)=L^p(X)$ of the set of all $X$-valued analytic polynomials:
$${\cal P}_+(X)=\left\{ \sum_{k=0}^{n} z^k x_k;n\in
\nat,z\in \tore \hbox{ \
and
\
}
x_k\in X\right\}. $$
\par
The corresponding spaces of scalar functions
 are denoted simply $L^p$
\
and
\
$H^p$.\par
\bg
These definitions  also make sense for 
\
$p\in ]0,1[ $ ,
\
but we only obtain quasi-Banach spaces instead of  Banach spaces .\par
\bg
If $p\geq 1$, 
then using the Poisson kernel we see easily that 
$\widetilde H^p(X)$
\
 can be identified with a subspce of
\
 $H^p(X)$,
but in general $\widetilde H^p(X) \not = H^p(X)$ .
 However,when $0< p< +\infty $ , equality occurs if
and only if $X$ has the analytic Radon-Nikodym property.
 For references and more information on 
this property
see [{\bf E1}], [{\bf E2}] and [{\bf GLM}].\par
\bg
If $X$ and $Y$ are two Banach spaces then ${\cal L}(X,Y) $
denotes the space of all bounded
 operators from
$X$ into $Y$. A norm one operator
$ u \in {\cal L}(X,Y)$
is $\lambda$-surjective, if for every
 $y \in Y $ of norm
smaller than one, there exists 
$x \in X$ of norm smaller than
$\lambda$
such that $u(x) = y$. A norm one operator is  said to be a metric
surjection, if it is $(1+ \eps)$-surjective for
every $\eps > 0$, \ie if its transpose is an
isometric embedding.\par
 \bg
If $T \in {\cal L}(X,Y)$ is an operator, then the formula
$$\widetilde T \left(\sum_{k=0}^n z^k a_k \right )=\sum_{k=0}^n z^k T(a_k)$$
defines an operator
 $\widetilde T \in {\cal L} (\widetilde H^p(X),\widetilde H^p(Y))$
of the same norm .\par
\bg
We recall the definition of the projective tensor norm $\N{.}_\wedge$ 
on 
$ X \otimes Y $.
$$ \forall u \in X \otimes Y ~~\N{u}_\wedge =\hbox{Inf}
\left\{\sum_{k=0}^n\N{x_k}\N{y_k}\right\} $$
where the
\
 infimum runs over all
\ 
possible representations 
\
of the form
\
$u =\sum_{k=0}^n x_k \otimes y_k$ ;
$x_k \in X $ 
and 
$y_k \in Y$.\par
\bg
We denote by $X \hat\otimes Y$
the completion  of  $X \otimes Y $ 
equipped with the projective norm .
The resulting Banach  space is called the projective tensor product
of $X$ and $Y$ (\cf[{\bf G1}],[{\bf G2}]).\par
\bg
 Let $J_{X,Y}$ be the natural linear operator of norm one from
$ H^2(X)\hat\otimes H^2(Y)$ 
into
$H^1(X\hat\otimes Y )$ defined by
$$ J_{X,Y}(f\otimes g)(z) = f(z)\otimes g(z) $$
we also  consider  the induced opertor $j_{X,Y}$
from  $\widetilde H^2(X)\hat\otimes \widetilde
H^2(Y)$ into $\widetilde H^1(X\hat\otimes Y )$ . \par
\bg
It is almost 
clear that if
$X$ 
or 
$Y$ is an
 $L^1$-space, then
$J_{X,Y}$ (resp.$j_{X,Y}$)
 is onto,
\ie $\lambda$-surjective 
for
some $\lambda$;
 this is a simple
 consequence of the scalar case.
\
On the other hand ,
if $X$
 and
 $Y$
 are Hilbert spaces then 
$J_{X,Y}$ 
is a metric surjection; 
\
this is a classical result due to
 Wiener-Masani
 in the finite-dimensional case, 
and to Sarason 
in the infinite-dimensional case. 
Recently G.Pisier in [{\bf P3}]
has generalized  the preceding results by proving the surjectivity
of
 $j_{X,Y}$,
 if the considered spaces have type 2,
or if they are
2-convex Banach lattices. 
For a more detailed discussion see [{\bf P3}].\par
\bg
The urgent question was then to find Banach spaces $X$ and $Y$
such that $j_{X,Y}$ is not onto.
\
This was answered in our earlier paper [{\bf K}] ,
where  $H^1$-projective 
Banach spaces (see definition below) 
 were used in an essential
way. Indeed ,the role played by  $H^1$-projective 
\
Banach spaces can be clarified by 
the following proposition from [{\bf K}].
\par
\bg 
\proclaim Proposition. Let $X$ and $Y$ be $H^1$-projective 
Banach spaces. Then the following assertions are
equivalent:
 \item{$i.$}\quad $X \hat\otimes Y $  is an $ 
H^1$-projective Banach space.
\item{$ii.$}\quad $j_{X,Y} $ is onto.
\par
\bg
This suggested to us to make 
a somehow 
systematic study of 
$H^1$-projective 
 Banach spaces, and this is the purpose of
this paper.
Before describing
the organization of the paper let us
recall some facts which will be
frequently used in the sequel.
\par
The following result comes from
[{\bf HP}].
\proclaim Proposition 1.1. If 
$\sigma \in {\cal L}(X,Y)$
is a (metric) surjection
then the following assertions
are equivalent.\sn
\item{$i.$}
There exists $p\in [1,+\infty]$  
such that the operator $
\tilde\sigma :
\widetilde H^p(X) \rightarrow
\widetilde H^p(Y) $
 is a 
  (metric) surjection. \sn
\item{$ii.$}
For all $ p\in ]0,+\infty]$,
  the operator $
~\tilde\sigma :
\widetilde H^p(X) \rightarrow
\widetilde H^p(Y) $
 is a (metric) surjection.
\par
We will  also need  the principle of
local reflexivity [{\bf LR}].
The following formulation comes
from[{\bf D}].
\par
\proclaim Proposition 1.2. Let $X$
be a Banach space, and 
let $F$ and $G$ 
be finite-dimensional subspaces
of $X^{\ast}$ and $X^{\ast\ast}$
respectively  then for every
$\eps > 0 $ there exists an isomorphism
$T : G \rightarrow T(G) \subset X$
with the following properties :\sn
\item{$i.$} For every $
x\in X\cap G  $  we have $  T(x)=x $.\sn
\item{$ii.$}
$\N{T}.\N{T^{-1}} \leq 1+\eps $.\sn
\item{$iii.$}  For every 
 $ y\in F $
 and every $  x\in G $  we have 
$\langle y,x\rangle =\langle
y,T(x)\rangle.$\sn
\par
Let $K_n $ denote the Fejer kernel
defined by
$$ K_n(\theta)={1\over n}\left ({\sin
n\theta /2}\over {\sin \theta /2}\right
)^2 .$$
We recall that the De La
Vall\'ee-Poussin kernel $V_{n,r}$
is defined by
$$ V_{n,r}={{rK_{rn}-K_n}\over {r-1}};~~
r=2,3,...$$
It  is  easy to see that $\widehat V_{n,r}(k)=1$ 
if
$\abs{k}\leq n$ and
$\widehat V_{n,r}(k)=0$
if
$\abs{k}\geq rn$,
 one also  has $\N{V_{n,r}}_1\leq 
{{r+1}\over {r-1}} $.
The convolution operator
$f \rightarrow f\ast V_{n,r}$
will be denoted $\Phi_{n,r} $.
See [{\bf Z}] for more information
on this topic .
\par
For any other Banach space terminology,
we refer the reader to any general
treatise on Banach spaces like [{\bf
LT}].
\par
Let us now describe the organization
 of this paper.
\par
In section 2. we give several 
equivalent
formulations of the
$ H^1$-projective  property.
\par
In section 3. we relate the
$H^1$-projective property 
to martingale inequalities and
prove that in an $H^1$-projective 
space $X$, Hardy martingales valued
in $X$ converge unconditionally
and that such a space has the 
analytic Radon-Nikodym property.
\par
In section 4. we prove that 
$H^1$-projective spaces have cotype 2
and satisfy Grothendieck's theorem.
\par
In section 5. we show that the 
ultraproduct of $H^1$-projective 
spaces is also $H^1$-projective.
\par
Finally in section 6. we give some
examples of $H^1$-projective spaces.\par
\bg
\noindent{\bf 2. H${}^{\fivebf 1}$-Projectivity,
equivalent formulations}
\bg
If $F$ is an element of  the  tensor
product $H^1 \otimes X$, then for every 
representation $\sum_0^n h_k \otimes
x_k$ of $F$ we have
$$\N{F}_{\widetilde H^1(X)}
\leq \sum_{k=0}^n \N{h_k}_{H^1} \N{x_k}_X $$
hence for every $F$ in  $H^1 \otimes X$
the following inequality holds
$$ \N{F}_{\widetilde H^1(X)}
\leq \N{F}_{ H^1\hat\otimes X}.$$
So we always have a norm one 
inclusion $H^1\hat\otimes X
\hookrightarrow \widetilde H^1(X) $ .
\par
Let us then make the following
definition:
\par
\noindent{\bf Definition 2.1}. A complex Banach space
$X$ is called $H^1$-projective if
the canonical inclusion  map $H^1
\hat\otimes X
\hookrightarrow \widetilde H^1(X) $ is
surjective.
\par
Equivalently, $X$ is $H^1$-projective
 if there exists some constant $c$
such that 
$$ \forall F \in H^1\hat\otimes X ,\qquad\N{F}_{H^1\hat\otimes X}
\leq c
\N{F}_{\widetilde H^1(X)} \eqno (2.1) $$
The smallest constant $c$ satisfying
(2.1) is called the $H^1$-projectivity
constant of $X$, and will be 
denoted $\eta (X) $.
\par
This property was already introduced
in our earlier paper [{\bf K}].
 In what follows we give some
characterizations of $H^1$-projective
spaces and $H^1$-projective dual
spaces.
\par
The following proposition
characterizes $H^1$-projectivity
in termes of the possibility of lifting
$X$-valued analytic functions. More
precisely we have:\par

\proclaim Proposition 2.2. Let
$X$ be a complex Banach space. The 
following assertions are equivalent :
\item{$i.$} $X$ is $H^1$-projective.
\item{$ii.$} There exist a constant $\lambda$ and 
a metric surjection  $\sigma : \ell^1(I) \rightarrow X$ such that $\tilde\sigma : \widetilde H^1(\ell^1(I))
\rightarrow \widetilde H^1(X)$ is a $\lambda$-surjection.
\item{$iii.$} Every surjection $\rho$  
from some Banach
space $Y$ onto $X $ induces a surjection 
$\tilde\rho : \widetilde H^1(Y)\rightarrow \widetilde H^1(X)$.
\item{} Moreover, the smallest $\lambda $ satisfying  $ii.$ is
equal to $\eta (X)$.\par

\noindent Proof : Note that if $S_X$ is the 
unit sphere of $X$ 
then the operator 
$\sigma \in {\cal L}(\ell^1(S_X),X)$
defined by
$$\sigma ((\alpha _x)_{x \in S_X})
=
\sum_{x \in S_X } \alpha _x . x $$
is clearly a metric surjection. Now
since we always have 
$\ell^1(I) \hat\otimes H^1
=
\ell^1(I,H^1) = \widetilde H^1(\ell^1(I)) $,  it is immediate 
that Im
$ \tilde\sigma = \widetilde H^1\hat\otimes
X$.
\par
$i.~~ \Leftrightarrow ii.$
This is obvious using the preceding
observation.
\par
$iii.\Rightarrow ii.$
Clear.
\par
$ii.~
\Rightarrow
iii.$ Using the 
lifting property of $\ell^1(I)$, 
if $\rho : Y\rightarrow X $ is any
$K$-surjection, we can find
$v : \ell^1(I)\rightarrow Y $
such that $\rho\rond v = \sigma$
and $\N{v} \leq K $.
\par
Let $f \in \widetilde H^1(X) .~$
Since $\tilde \sigma$
is a $\lambda$-surjection, we can
find $g \in \widetilde H^1(\ell^1(I))$
with $\tilde\sigma(g)=f$
and $\N{g}\leq \lambda\N{f}$.
Consider then $h=\tilde v(g) \in \tilde
H^1(Y)$ obviously 
$\tilde\rho (h)=f$ and 
$\N{h}\leq \lambda K \N{f}$.
\cqfd\bg
\noindent{\it Remark.}
 We can replace in the preceding
proposition $\widetilde H^1$
by 
$\widetilde H^p$
for any $p \in [1,+\infty] $, this
follows from proposition 1.1.
\bg
Proposition 2.3, below, expresses
$H^1$-projectivity
by means of 
the possibility 
of extending
some
operators defined on 
$H^1$. 
To this 
end we  need the 
following notation.
We will say that an operator 
$u \in {\cal L}(L^1,X)$
(resp.
${\cal L}(H^1,X)$)
belongs to 
$ {\cal L}_F(L^1,X)$
(resp.
${\cal L}_F(H^1,X)$)
if there exists a positive integer
$k_u$ such that
$u(e^{in(.)}) = 0 $ for every 
$n\not\in [-k_u,k_u]$
(resp. $n>k_u $).
\par

\proclaim Proposition 2.3. For every
complex Banach space $X$, the following
assertions are equivalent :
\item{$ i.$} $X $ is $H^1$-projective.
\item{$ii.$} There exists $ \lambda > 0 $  such
that every  $v\in {\cal L}(H^1,X^{\ast})$  
has an extention $\bar v \in
 {\cal L}(L^1,X^{\ast})$, \ie  $\bar v_{\vert
H^1} =v$, with $\N{\bar v}\leq \lambda \N{v}$. 
\item{$iii.$} There exists $\lambda > 0 $  such
that every $ v\in
 {\cal L}_F(H^1,X^{\ast}) $ 
has an extention $ \bar v \in 
 {\cal L}_F(L^1,X^{\ast}) $, with $
\N{\bar v}\leq \lambda \N{v} $.
\item{}Moreover, the smallest constant
satisfying $ii$ (resp. $iii$) is equal to $\eta(X)$.
\par

\noindent Proof : $i . \Leftrightarrow 
ii .$
This is just duality.
Indeed, 
letting
$I:H^1\hat\otimes X
\rightarrow
\widetilde H^1(X) $
be the canonical inclusion map,
$H^1$-projectivity
of $X$ 
is by definition
the statement that $~I~$ 
is an isomorphism or 
equivalently that 
$I^{\ast}$
is an isomorphism.
But
$$ (H^1 \hat\otimes X)^{\ast}
=
{\cal L}(H^1,X^{\ast}) .$$
And
$$(\widetilde H^1(X))^{\ast}
=
{\cal L}(L^1,X^{\ast})/
[\widetilde H^1(X)]^{\perp} $$
with
$[\widetilde H^1(X)]^{\perp}
=\{ v\in {\cal L}(L^1,X^{\ast})~:~v_{\vert H^1} = 0 \}$.
\par
So $H^1$-projectivity
of $X$ is equivalent to
the statement
that 
${\cal R}:
{\cal L}(L^1,X^{\ast})
\rightarrow
{\cal L}(H^1,X^{\ast})
: ~v\mapsto 
v_{\vert H^1}$
is onto and this is $ii .$
\par
$ii .
\Rightarrow
iii .
\Rightarrow
i.$
This is easy
using the operators
$\Phi_{n,r}$
(see section 1).
Details are left as an exercise
to the reader.\cqfd
\par
In the next theorem we show that
$X$
is $H^1$-projective
if and only if 
every operator in
${\cal L}(X^{\ast},H^{\infty})$
factors through some
$\ell^{\infty}(I)$-space.
\par
\proclaim Theorem 2.4. 
Let $X$ be a complex Banach
space. The following assertions 
are equivalent :
\item{$i.$} $X$ is $H^1$-projective.
\item{$ii.$}
There exist a constant $K$  and an
isometric embedding $j:X^{\ast} \rightarrow \ell^{\infty}(I)$ such 
that every operator $u \in {\cal L}(X^{\ast}, H^{\infty})$ extends to an operator $\bar u \in {\cal L}(\ell^{\infty}(I),
H^{\infty}) $ with $ \bar u\rond j= u$ and $
\N{\bar u}  \leq K \N{u}$.
\item{$iii.$}
 There exists a constant 
 $K$   such that for every
 isometric  embedding $l:X^{\ast} \longrightarrow  Z $
 of  $X^{\ast}$ 
 into a Banach space  $Z$; 
 every operator 
$u \in {\cal L}(X^{\ast}, H^{\infty})$ 
 extends to an operator 
$ \bar u \in {\cal L}(Z, H^{\infty}) $
  with 
 $ \bar u\rond l= u $
  and 
 $ \N{\bar u}  \leq K \N{u} $.
\par
\noindent Proof : $i\Rightarrow ii$.
By Proposition 2.2 and 
Proposition 1.1
we can  find a metric surjection
$\sigma : \ell^1(I) \rightarrow X$
such that 
$\tilde\sigma : \widetilde H^{\infty}
(\ell^1(I))
\rightarrow \widetilde H^{\infty}(X) $
is a $K_1$-surjection for some
$K_1$.
Cleary $j=\sigma^{\ast} :
X^{\ast}\rightarrow
\ell^{\infty}(I)$
is an isometric embedding.
\par
Consider $u \in {\cal L}(X^{\ast},
H^{\infty})$
and  put 
$u_n=\Phi_{n,2}\rond u$
where
$\Phi_{n,2}:
H^{\infty}\rightarrow
H^{\infty}:
f \mapsto f\ast V_{n,2}$.
Obviously
$u_n$
takes its
values in the span of
$\{ e^{ik(.)}\}_{0\leq k\leq 2n} $
so we can find,
for each $n$,
a sequence
$x_0^{(n)},...,x_{2n}^{(n)}$
such that
$$\forall x^{\ast} \in X^{\ast},\qquad
 u_n(x^{\ast})=
\sum_{k=0}^{2n} x^{\ast}
(x_k^{(n)})e^{ik(.)} $$
So if 
$f_n=\sum_{k=0}^{2n}
e^{ik(.)}x_k^{(n)} ~
\in \widetilde H^{\infty}(X)$
we obtain
$$\N{f_n}_{\widetilde H^{\infty}}
=\N{u_n} \leq \N{u} \N{V_{n,2}}
\leq 3 \N{u} \eqno (2.2) $$
and
$$\forall x^{\ast} \in X^{\ast},\qquad
 u_n(x^{\ast})=
\tilde x^{\ast}(f_n) . \eqno (2.3)$$
Using the hypothesis we can
find $g_n \in \widetilde H^{\infty}
(\ell^{\infty}(I))$
such that 
$$\tilde\sigma (g_n)=f_n
\hbox{
and
}
\N{g_n} \leq K_1 \N{f_n}. \eqno (2.4)$$
Define 
$\bar u_n \in {\cal L}
(\ell^{\infty}(I),H^{\infty})$
by
$\bar u_n(z^{\ast})=\tilde z^{\ast} (g_n)$
 for every  $z^{\ast} \in (\ell^1(I))^{\ast}= \ell^{\infty}(I)$.
Obviously by (2.4) and (2.2) we obtain
$$\forall n,\qquad
\bar u_n \rond j = u_n \quad
\hbox{and}\quad
\N{\bar u_n} \leq 3 K_1 \N{u}.
\eqno  (2.5)$$
Let now ${\cal U}$ be a non 
trivial ultrafilter on $\nat$.
Using the fact that
$(L^1/H^1)^{\ast} = H^{\infty}$,
we can define for each $t\in
\ell^{\infty}(I) $
the following limit in the
weak$-\ast$
topology
$\bar u(t) =\lim_{\cal U}
\bar u_n(t) $.
We conclude immediately
that $\bar u \in {\cal L}(\ell^{\infty}(I),
H^{\infty})$
with
$\N{\bar u}  \leq 3 K_1 \N{u}$.
On the other hand, for each 
$x^{\ast} \in X^{\ast}$
and each $f \in L^1$ 
we have
$$\eqalign{\langle\bar u\rond j(
x^{\ast}),f\rangle &=
\lim_{\cal U}~\langle \bar u_n\rond j(
x^{\ast}), f\rangle \cr
&=\lim_{\cal U}~\langle u_n(
x^{\ast}),f\rangle \cr
&=\lim_{\cal U}~\langle V_{n,2}\ast u(
x^{\ast}),f\rangle \cr
&=\lim_{\cal U}~\langle u(
x^{\ast}),V_{n,2}\ast f\rangle \cr
&=\langle u(
x^{\ast}),f\rangle .} $$
The last equality holds since
$\{f\ast V_{n,2} \}_n $
converges to $ f $ in the $L^1$-norm.
We conclude that
$\bar u \rond j = u $
and $ii.$ follows.
\par
$ii.
\Rightarrow
iii.$
Let $l:X^{\ast}\rightarrow Z$
be an isometric embedding, then
by the Hahn-Banach theorem
we can find $\bar \jmath:Z \rightarrow 
\ell^{\infty}(I) $
such that $\bar \jmath\rond l=j$
and $\N{\bar \jmath}=\N{j}=1 $.
\par
Consider $u \in {\cal L}(X^{\ast},
H^{\infty})$,
there is $\bar u_1 \in {\cal L}(
\ell^{\infty}(I),
H^{\infty})$
such that
$\bar u_1 \rond j = u$
and $\N{\bar u_1} \leq K \N{u}$.
Define then 
$\bar u=\bar u_1\rond \bar \jmath
 \in {\cal L}(Z,
H^{\infty})$.
Clearly $\bar u \rond l = u$
and
$\N{\bar u } \leq K \N{u} $.
\par
$iii.\Rightarrow i.$
Let $\sigma : \ell^1(I)
\rightarrow  X$
be any metric surjection. 
We will show that
$\tilde\sigma : \widetilde H^{\infty}
(\ell^1(I))
\rightarrow \widetilde H^{\infty}(X) $
is a surjection.
This implies the  result,
by Propositions 1.1 and 2.2 .
\par
Let $f\in \widetilde H^{\infty}(X)$,
and consider the operator 
$u: X^{\ast}\rightarrow H^{\infty}$
defined by
$u(x^{\ast})=\tilde x^{\ast}(f)$.
Clearly $\N{u} = \N{f}_{\tilde
H^{\infty}(X)}$
and by $iii.$ we
can find
$\bar u :\ell^{\infty}(I)
\rightarrow H^{\infty} $
such that
$$\bar u\rond \sigma^{\ast}= u
\quad\hbox{and}\quad
\N{\bar u}\leq K\N{u}=K\N{f}_{\tilde
H^{\infty}(X)}.\eqno (2.6) $$
Let $e_{i} \in \ell^{\infty}(I)$
be defined by
$e_{i}(j ) =0$
if
$j \not = i$
and
$e_{i}(j) =1$
if
$j = i$,
and put $g_{i}=\bar u(e_{i})
\in H^{\infty}$.
It is easy to see that 
$g=(g_{i})_{i \in I}
\in \widetilde H^{\infty}(\ell^1(I))$
and that
$$
\forall z^{\ast}\in (\ell^1(I))^{\ast},
\qquad\bar u(z^{\ast})=\tilde z^{\ast}
(g)$$
moreover,
$\N{g}_{\widetilde H^{\infty}(\ell^1(I))}
= \N{\bar u} \leq K \N{f}_{\tilde
H^{\infty}(X)} $
by (2.6).
On the other hand, the equality
$\bar u\rond \sigma^{\ast}= u$
from (2.6) is equivalent to 
$\tilde \sigma (g) = f$.
This completes the proof of the theorem.\cqfd
\par
\bg
In what follows we give a
characterization
of $H^1$-projective dual spaces.
\par
\proclaim Proposition 2.5. Let $X$ be a complex Banach
space. The following assertions are
equivalent :
\item{$i.$} $X^{\ast}$  is $ H^1$-projective.
\item{$ii.$} There exists a  constant $\lambda$
  such that every  operator  $ v \in {\cal L}_F (H^1,X) $ has an extention 
$\bar v \in {\cal L}_F (L^1,X) $  such that $\N{\bar v} \leq \lambda \N{v} $.
\item{}Moreover, the smallest $\lambda$ satisfying $ii.$
is equal to $\eta (X^{\ast})$.
\par
\noindent Proof : $i.\Rightarrow ii.$
Let
$j:X \rightarrow X^{\ast}$
be the natural
embedding.
Consider  
$v \in {\cal L}_F (H^1,X)$
and $\eps > 0$.
\par
By proposition 2.3.
$iii$
there exists an extention
$v_1 \in {\cal L}_F
(L^1,X^{\ast\ast})$
of
$j\rond v$
such that
$\N{v_1} \leq \lambda \N{v} $.
\par
Let $G~=$ Im $~v_1 \subset
X^{\ast\ast}$
and
$F=\{0\}$;
using the local reflexivity
principle (see Proposition 1.2),
we find
$T:G\rightarrow X$ such that
$\N{T} \leq 1+\eps$
\
and
$T_{\vert G\cap X}=
id_{\vert G\cap X}$.
So if we put
$\bar v=T\rond v_1$,
we obtain
$\bar v \in {\cal L}_F
(L^1,X)$ such that
$\N{\bar v} \leq (1+\eps)\lambda
\N{v}$
and
$\bar v_{\vert H^1} =v$.
This last equality
 comes from the fact
$\hbox{Im}~v_{1\vert H^1} \subset
G\cap X$.
\par
$ii.
\Rightarrow
i.$
We need the following lemma.
\par
\proclaim Lemma. The set
$\{ v \in {\cal L}_F
(H^1,X)~:~\N{v} \leq 1\} $
is dense 
in the unit ball of
${\cal L}
(H^1,X^{\ast\ast})$
with respect to the
weak-$\ast$
topology
(\ie $\sigma ({\cal L}
(H^1,X^{\ast\ast}),
H^1\otimes X^{\ast})$).
\par
\noindent Proof of the lemma :
If the lemma is false, then by the
Hahn-Banach theorem,
we can find $u$ in the closed unit
ball of 
${\cal L}
(H^1,X^{\ast\ast})$, 
 $f $ in $ H^1\hat\otimes X^{\ast}$,
and
$\eps >0$
such that
$$\forall v \in {\cal L}_F(H^1,X ),~
\N{v}\leq 1 \Rightarrow
\hbox{Re}(\langle u,f\rangle -
\langle v,f\rangle) \geq 4\eps
\eqno (2.7)$$
Since the set of analytic
$X^{\ast}$-valued polynomials
${\cal P}_+(X^{\ast}) $
is norm-dense   in
$H^1\hat\otimes X^{\ast} $, 
we can find
$g=\sum_0^n  e^{ik(.)}~ b_k
\in {\cal P}_+(X^{\ast}) $
such that (2.7)
becomes
$$\forall v \in {\cal L}_F(H^1,X ),~
\N{v}\leq 1 \Rightarrow
\hbox{Re}(\langle u,g\rangle -
\langle v,g\rangle) \geq 2\eps
\eqno (2.8)$$
Define $ r>1$ by the condition
$${3\over {r + 2}}\abs{
\langle u,g\rangle} < \eps \eqno
(2.9)$$
and consider
$u_1 = {{r-1}\over {r+1}} u\rond
\Phi_{n,r}$
(see section 1 for the notation),
clearly $u_1\in 
{\cal L}_F(H^1,X^{\ast\ast} )$
and
$\N{u_1}\leq {{r+1}\over {r+2}}
\N{u} <1$.
Moreover,
using (2.8) and (2.9)
we have for every 
$v \in {\cal L}_F(H^1,X )$
of norm smaller than one:
$$\hbox{Re}(\langle u_1,g\rangle -
\langle v,g\rangle) \geq 
\hbox{Re}(\langle u,g\rangle -
\langle v,g\rangle) 
- {3\over {r+2}}\abs{\langle
u,g\rangle}> \eps .$$
So we have found  
$u_1$
 in the open unit ball of 
${\cal L}_F(H^1,X^{\ast\ast} )$
and
$g=\sum_0^n  e^{ik(.)}~ b_k
\in {\cal P}_+(X^{\ast}) $
such that
$$\forall v \in {\cal L}_F(H^1,X ),~
\N{v}\leq 1 \Rightarrow
\hbox{Re}(\langle u_1,g\rangle -
\langle v,g\rangle) \geq \eps.
\eqno (2.10)$$
Consider now the finite-dimentional
spaces
$$F=\hbox{Span}\{ b_k :~0\leq k\leq n
\} \subset X^{\ast}\quad\hbox{
and}\quad
G=\hbox{Span}\{u_1(e^{ik(.)}) :~
k\geq 0 \} \subset X^{\ast\ast}$$
By the local reflexivity
principle there
exists an operator 
$T: G\rightarrow X$
 satisfying\par
\item{$i.$} \quad$T_{\vert G\cap X}=id_{\vert G\cap X}$.
\item{$ii.$} \quad$\N{T} \leq {1\over {\N{u_1}}}$.
\item{$iii.$} \quad$\forall y\in F,\quad\forall z \in G,~\langle z,y\rangle =
\langle T(z),y\rangle$.\par
We define
then  $v \in {\cal L}_F(H^1,X ) $
by
$v=T\rond u_1$.
It is easily seen that
$\N{v} \leq 1$ 
and
$\langle u_1,g\rangle =
\langle v,g\rangle $.
This contradicts (2.10)
and proves the lemma.\cqfd
\par

We can now finish the proof of
$ii.
\Rightarrow i.$
Take
$f=\sum_0^n  e^{ik(.)}~ a_k$
an element of
$H^1\otimes X^{\ast}$,
then the following holds
$$\eqalign{\N{f}_{H^1\hat\otimes 
X^{\ast}}&=
\hbox{ Sup }\{\abs{\langle u,f\rangle}
:~u\in{\cal L}
(H^1,X^{\ast\ast}) ,~\N{u} \leq 1\}\cr
&=\hbox{ Sup }\{\abs{\langle u,f\rangle}
:~u\in{\cal L}_F
(H^1,X) ,~\N{u} \leq 1\}\cr
&\leq \lambda
\N{f}_{\widetilde H^1(X^{\ast})},}$$
the second equality comes from
the lemma,
and the last inequality comes from
the hypothesis.
This, of course, achieves the proof.\cqfd
\par
\bg
\proclaim Corollary 2.6. A complex
Banach space $X$ is $H^1$-projective, 
if and only if, 
$X^{\ast\ast}$
is $H^1$-projective .
Moreover,  $\eta (X) = \eta
(X^{\ast\ast})$.
\par
This is a direct consequence
of Proposition 2.3 and
Proposition 2.5.
\par\bg
\noindent{\it
Remark.}
Let $j:X\rightarrow \ell^{\infty}(J)$
be any isometric
embedding, then $\sigma = j^{\ast} :
(\ell^{\infty}(J))^{\ast}
\rightarrow
X^{\ast}$
is a metric surjection.
Since
 $(\ell^{\infty}(J))^{\ast}$
is 
an abstract
$L^1$-space, 
it is isometric to some $L^1(\mu)$.
Now if $X^{\ast}$ 
is $H^1$-projective,
then
$\tilde\sigma : \widetilde H^{\infty}
(L^1(\mu))\rightarrow
\widetilde H^{\infty}(X^{\ast}) $
is a surjection.  
 But using the fact that
$\sigma = j^{\ast} $
and a remark from [{\bf HP}]
we see that 
$\tilde\sigma : H^{\infty}
(L^1(\mu))\rightarrow
 H^{\infty}(X^{\ast}) $
is also a surjection. So we obtain 
the follwing :
\bg
$X^{\ast}$ is $H^1$-projective, 
if and only if, there exists a
metric surjection
$\sigma : L^1(\mu)\rightarrow
X^{\ast}$
such that 
$\tilde\sigma : H^{\infty}
(L^1(\mu))\rightarrow
H^{\infty}(X^{\ast}) $ 
is a surjection.

Note that this does not follow
from corollary 3.2 below.

\bg
\noindent{\bf 3. Martingale
inequalities in $H^1$-projective
spaces } \par
\bg

Let us start this section
 by recalling some definitions
and notation about 
Hardy martingales.
The notion of Hardy martingales
appeared first in
[{\bf Ga}], 
but was already implicit in [{\bf
E2}].
\par
We consider the infinite dimensional
torus 
$\tore^\nat$
and denote by $\theta_n$ the
$n\hbox{th}$ coordinate of a
point
 $\Theta$
in $\tore^\nat$.
Let
$(\Omega,{\cal A},\pe)$
denote $\tore^\nat$
equipped with its
normalized Haar measure.
We denote by ${\cal A}_n$
the $\sigma$-algebra generated
by $(\theta_1,...,\theta_n)$
on $\tore^\nat$.
Let $(M_n)_{n\geq 0}$
be a sequence in 
$L^1(\Omega,{\cal A},\pe;X) $
which is a martingale with respect
to
the sequence of $\sigma$-algebras
$({\cal A}_n)_{n\geq 0} $. 
We will say that 
$(M_n)_{n\geq 0}$
is a Hardy martingale,
if for each fixed
$\theta_1,\theta_2,...,\theta_{n-1}$
the function
$$\theta \mapsto M_n(
\theta_1,...,\theta_{n-1},\theta)
\eqno (3.1) $$
is in $\widetilde H^1(X)$.
Let $\hbox{d}M_n= M_n-M_{n-1}$
for $n \geq 1$,and
$\hbox{d}M_0= M_0$.
Equivalently we require that
$\hbox{d}M_n$
satisfies
$$\forall k\geq 1 ~~~
\int_\tore \hbox{d}M_n(
\theta_1,...,\theta_{n-1},\theta)~
e^{ik\theta} dm(\theta) =0. $$
Finally, if the function
in (3.1) is always of the form
$x+e^{i\theta}y$ 
for some $x,y \in X$, 
then $(M_n)_{n\geq 0}$
is called an analytic martingale.
\par
The following proposition is a
combination of results from
[{\bf GM}] and [{\bf HP}].
\par

\proclaim Proposition 3.1.
Let $X$ be an $H^1$-projective 
space. Then every 
$L^1$-bounded
$X$-valued
Hardy martingale satisfies
$$\forall m\geq 1,\qquad
\left [\sum_{n=1}^m
(\ee\N{\hbox{d}M_n})^2\right
]^{1\over 2} \leq 2 \eta (X)
\ee \N{M_m} .\eqno (3.2) $$
\par
\noindent Proof :
It is sufficient to prove the
result for Hardy martingales
$(M_n)_{n\geq 0}$
whose differences 
$\hbox{d}M_n(
\theta_1,...,\theta_n)$
are $X$-valued
trigonometric polynomiales in 
the variables
 $\theta_1,...,\theta_n $.
So we assume this in the sequel.
\par
Let $\eps$ be an arbitrary
positive number,
and $m$ a positive integer .
By a result of [{\bf GM}]
there exist an analytic function
$F:{\cal D}\rightarrow X$,
a continuous function 
$p:\tore \rightarrow \tore^\nat $,
and an increasing sequence
$0<r_0<r_1<...<r_n <...<1 $
such that 
$$\eqalign { 1&.\quad p^\ast (dm( \theta))
=\Otimes_{n=1}^\infty
dm(\theta_n) \cr
2&.\quad F(0)= M_0 \cr
3&.\quad \forall n,\forall \theta \in\tore,
\qquad\N{F(r_{n+1}e^{i\theta})
-F(r_n e^{i\theta})
 - \hbox{d}M_{n+1}(p(\theta))}_X
\leq \eps }\eqno (3.3)$$
Define then $G(\theta)=F(r_m
e^{i\theta})$,
so that we have
$G\in \widetilde H^1(X) $. 
On the other hand, by Proposition 
2.2, there exists a metric surjection
$\sigma :\ell^1 (I) \rightarrow X$
such that 
$\tilde\sigma :\widetilde H^1(\ell^1(I))
\rightarrow
\widetilde H^1(X)$
is an $\eta (X)$-surjection . 
So we can find $\tilde G \in 
\widetilde H^1(\ell^1(I))$ 
satisfying
$$\tilde\sigma (\tilde G)=G 
\quad\hbox{
and
}\quad
\N{\tilde G}_{\widetilde H^1(\ell^1(I))}
\leq \eta (X)
\N{G}_{\widetilde H^1(X)} \eqno (3.4) $$
Let $P_r$ be the Poisson kernel,
and let $\tilde G_r $ 
denote $\tilde G \ast P_r $.
We know using a result from
[{\bf HP}] that
$$
\left (\sum_{n=0}^{m-1}
\N{\tilde G_{\rho_{n+1}} -
\tilde G_{\rho_n} }_{\tilde
H^1(\ell^1(I))}^2 \right )^{1\over 2}
\leq 2 \N{\tilde G_{\rho_m}}_
{\widetilde H^1(\ell^1(I))}
\leq 2 \N{\tilde G}_
{\widetilde H^1(\ell^1(I))}$$
for every sequence
$\rho_0<\rho_1<...<\rho_m \leq 1 $.
\par
Using (3.4), we can go back to $G$
and obtain 
$$
\left (\sum_{n=0}^{m-1}
\N{ G_{\rho_{n+1}} -
 G_{\rho_n} }_{\tilde
H^1(X)}^2 \right )^{1\over 2}
\leq 2 \eta (X) \N{ G}_
{\widetilde H^1(X)}. \eqno (3.5) $$
Taking $\rho = {{r_n}\over {r_m}}$
for $0 \leq n \leq m $,
we can write (3.5) as follows
$$
\left[ \sum_{n=0}^{m-1}
\left( \int_\tore
\N{F(r_{n+1}e^{i\theta})
-F(r_n e^{i\theta}) }_X dm(\theta)
\right )^2 \right ]^{1\over 2}
\leq 2\eta (X)\int_\tore
\N{F(r_m e^{i\theta}) }_X
 dm(\theta) $$
taking into account (3.3) we
conclude that\par
$$\eqalign{
\left [ \sum_{n=0}^{m-1}
\left ( \int_\tore
\abs{~\N{\hbox{d}M_{n+1}(p(\theta))~}_X
- \eps} \, dm(\theta) \right )^2
\right ]^{1\over 2} 
&\leq 2 \eta (X) 
 \int_\tore
(\N{M_m (p(\theta))~}_X + m \eps)
\, dm(\theta)  \cr
\left [ \sum_{n=0}^{m-1}
(\ee\N{\hbox{d}M_{n+1}}_X )^2
\right ]^{1\over 2}
&\leq 2 \eta (X) \ee \N{M_m}_X
+ (2\eta (X)m + \sqrt{m}~)\eps}$$
Since $\eps$ 
is arbitrary  we obtain (3.2) by 
letting $\eps \rightarrow 0 $.\cqfd
\bg
Following Xu in [{\bf X}], a
quasi-Banach space $X$ will be called
Hardy-convexifiable, if it has an
equivalent quasi-norm $\N{.}_X$
such that, for some $p > 0 $
and $c>0$,we have
$$\forall f\in {\cal P}_+(X),
\quad\N{f(0)}_X^p + c
\N{f - f(0)}_{H^1(X)}^p 
\leq \N{f}_{H^1(X)}^p. $$
It is shown in [{\bf X}] that if
a Banach space   $X$ satisfies 
a martingale inequality like (3.2),
then it is Hardy-convexifiable
and in particular it has 
the analytic Radon-Nikodym property.
So we have the following corollary.
\bg
\proclaim Corollary 3.2. 
Every $H^1$-projective Banach space
$X$
has an equivalent Hardy-convex
 quasi-norm. In particular $X$
has the super-analytic Radon-Nikodym
property.
\bg
Using  this corollary, one sees
that the $H^1$-projectivity of 
$X$ is equivalent to the 
apparently  stronger  property :
$ H^1\hat\otimes X = H^1(X) $.  
\bg
In order to state the next
proposition we need to the following
definition.
\bg
\noindent{\bf Definition 3.3} [{\bf Ga}].
We will  say that $X$ has the 
property U.H.M.D. (Unconditionality
of Hardy Martingale Differences), if
there exists a constant $K$  such
that for every  Hardy martingale
$(M_n)_{n\geq 0} $
and every choice of signs
$\eps \in \{+1,-1\}^\nat $, 
we have
$$
\forall m,\qquad
\ee \N{\sum _{n=0}^m \eps_n
\hbox{d}M_n }_X \leq K \ee
\N{M_m}_X  .$$
\par
\bg
\proclaim Theorem 3.4. 
Every $H^1$-projective 
Banach space $X$
has the property U.H.M.D.\par
\bg
\noindent Proof : 
The proof is decomposed into 
two steps which involve classical
arguments.
We first prove the following:
\par
{\it Step} 1. {\sl
Let $X$ be a complex 
 Banach space.
Suppose that  there exist two 
sequences  of integers
$(a_n)_{n \geq 0} $ 
and
$(\lambda _n)_{n \geq 0} $ 
 such that 
$\forall n~~a_n\lambda _n < a_{n+1},~
 \lim_{n\rightarrow \infty} 
\lambda _n = +\infty $;
and that there exists some constant
$K$, 
such that for every function
$f \in \widetilde H^1(X) $
satisfying 
$\hat f(n) = 0 , \forall n \not 
\in \cup_k [a_k, \lambda_k a_k]$  
and for every 
$\eps \in \{ +1,-1 \}^\nat$ , 
the following inequality holds
$$
\N{\sum_{n=0}^\infty \eps_n
\left ( \sum_{a_n \leq r \leq
\lambda_n a_n } \hat f(r) e^{ir(.)}
\right ) }_{\widetilde H^1(X) }
\leq K \N{f}_{\widetilde H^1(X) }.
\eqno (3.6) $$
Then X is U.M.H.D.}
\bg
\noindent Proof of step 1:
 Consider the Hardy martingale
$(M_n)_{n\geq 0}$ defind by
$$
\eqalign {\hbox{d}M_0 &= M_0 ,\cr
\hbox{d}M_n &=
\sum_{(p_1,..,p_n) \in A_n} 
 M_{n,p_1,..,p_n}~ 
\hbox{exp}\left (i\sum_{k=1}^n p_k
\theta_k \right ) } \eqno (3.7) $$
where $A_n= \{ (p_1,..,p_n)
\in \ent^n :~
p_n > 0~\hbox{ and } ~ \forall k\leq n
~\abs{p_k} \leq r_n \}$. 
It is clearly sufficient to prove
the result for such martingales. 
\par
We construct 
inductively two
sequences of integers
$(\mu_n)_{n\geq 1},(\ell_n)_{n\geq
1} ~$ satisfying:
\par
\itemitem{$1.$}\quad $ (\ell_n)_{n\geq
1}$  is strictly increasing.
\par
\itemitem{$2.$}\quad $\forall~n,\forall~(p_1,...,p_n)
\in A_n$  we have
$\ds
a_{\ell_n} \leq \sum_{k=1}^n
\mu_k p_k \leq \lambda_{\ell_n}
a_{\ell_n}.$\hfill$(3.8)$\par
To start the induction we 
find $\ell_1$ 
satisfying $ \lambda_{\ell_1} > r_1$ 
and put $ \mu_1 = a_{\ell_1} $, 
so that $a_{\ell_1} \leq \mu_1
\leq r_1 \mu_1 < \lambda_{\ell_1}
\mu_1 $. 
\par
Suppose
$\ell_1,...,\ell_{n-1};\mu_1,...,
\mu_{n-1} $ 
are already constructed.
Since $\lim_{n\rightarrow \infty}
\lambda_n = + \infty $, 
 we can find $\ell_n > \max ~(r_n, \ell_{n-1}) $ 
such that 
$$\left (
{{\lambda_{\ell_n}}\over
{r_n}} - 1 \right ) a_{\ell_n}
\geq 2 + ( r_n +1) \sum_{k=1}^{n-1}
\mu_k .$$
This yields that
$$
{{\lambda_{\ell_n}}\over {r_n}}
a_{\ell_n} -\sum_{k=1}^{n-1}
\mu_k  \geq 2 + a_{\ell_n} +
r_n  \sum_{k=1}^{n-1}
\mu_k .$$
So we can find $\mu_n $ such that
$$
\mu_n \in [ a_{\ell_n} +
r_n  \sum_{k=1}^{n-1} \mu_k ,
{{\lambda_{\ell_n}}\over {r_n}}
a_{\ell_n} -\sum_{k=1}^{n-1}
\mu_k ] \cap \nat .$$
Equivalently,
$$
 a_{\ell_n} +
r_n  \sum_{k=1}^{n-1} \mu_k 
\leq \mu_n \leq r_n \mu_n \leq
 \lambda_{\ell_n} a_{\ell_n} -
r_n  \sum_{k=1}^{n-1} \mu_k .$$
This gives (3.8), and achieves the
construction.
\par
Now, fix $\Theta = (\theta_n)_{n\geq
1} \in \tore^\nat $ 
and consider $ \eta \mapsto
G_m(\Theta,\eta) ,~ \eta \in \tore$ 
defined by 
$$
G_m(\Theta,\eta) = M_m(
\theta_1+\mu_1
\eta,...,\theta_n+\mu_n \eta ) .$$
Using (3.8), it is easily seen that
$G_m(\Theta,.)$ 
is an element of $\widetilde H^1(X) $ 
satisfying the hypothesis of step 1.
So for every $\eps \in
\{+1,-1\}^\nat$,  
we have
$$
\int_\tore\N{\sum_{n=0}^m \eps_n
\left (\sum_{A_n}
 M_{n,p_1,..,p_n}~ 
\hbox{exp}\left (i\sum_{k=1}^n p_k
(\theta_k + \mu_k \eta)\right
)\right ) }_X\,dm(\eta) \leq K \int_\tore\N{
G_m(\Theta,\eta)}_X \,dm(\eta)$$
It is then sufficient to integrate
with respect to $\Theta$ and to 
use the translation invariance
of the Haar measure on $\tore^\nat$ 
 to obtain
$$
\ee \N{\sum _{n=0}^m \eps_n
\hbox{d}M_n }_X \leq K \ee
\N{M_m}_X . $$
\par
{\it Step} 2. {\sl
For every $H^1$-projective 
space $X$  and every sequence
$(\lambda_n)_{n\geq 0} $  
such that $\lim_{n\rightarrow\infty}
\lambda_n = +\infty $, 
there exist a sequence  
$(a_n)_{n\geq 0} $  
such that $\forall n,\quad a_{n+1} >
\lambda_n a_n $,  
and a constant $K$, with the
following property:    
 for every function
$f \in \widetilde H^1(X) $
satisfying 
$\hat f(n) = 0$,  for all $ n \not 
\in \cup_k [a_k, \lambda_k a_k]$  
and for every 
$\eps \in \{ +1,-1 \}^\nat$, 
we have
$$
\N{\sum_{n=0}^\infty \eps_n
\left ( \sum_{a_n \leq r \leq
\lambda_n a_n } \hat f(r) e^{ir(.)}
\right ) }_{\widetilde H^1(X) }
\leq K \eta (X) \N{f}_{\widetilde H^1(X)
}.$$
}\par
This step is proved just as in the
classical case of $X= \comp $.  
See for instance [{\bf CW}.Theorem
2.1].\cqfd
\bg
\noindent{\bf 4. H${}^1$-projectivity,
Grothendieck's theorem and cotype 2
spaces}
\bg
Let us first recall some basic
  definitions. We refer the reader to [{\bf P1}]
  for more  details  and  other references. 
\par
Let $G = \{+1,-1\}^\nat$, 
let $\mu$
be the uniform probability measure
on $ G$, and let  $ 
\eps_n : G \rightarrow \{+1,-1\} $    
be the $n$\up{th}  coordinate . 
A Banach  space  $X$  is  said
to be of cotype 2, if there is  a
constant $c $  such that for  all
$ x_1,...,x_n $ in $X$ we have  
$$
\left ( \sum_{k=1}^n \N{x_k}_X^2
\right )^{1/2} \leq  c 
\left ( \ee \N{\sum_{k=1}^n
\eps_k x_k }_X^2  \right )^{1/2
}.$$
We  denote  $c_2(X)$ the smallest 
constant  $c$ for which this holds.
\par
\noindent We say that an operator
$ u \in {\cal L}(X,Y) $ is  
1-summing, if there is a constant $c$
  such that for all $
x_1,...,x_n $ in $X$ we have 
$$
\sum_{k=1}^n \N{u(x_k)}_Y \leq c~ 
\hbox{Sup}~\left \{ \sum_{k=1}^n
 \abs{\langle \xi,x_k\rangle} :
~\xi \in X^\ast, ~\N{\xi} \leq
1 \right\}.$$
\noindent We denote $ \pi_1(u)$ 
the smallest constant $c$ with
  this property.
\par
We can now give the following 
definition.
\par
\noindent{\bf Definition 4.1} [{\bf P1}].  
We will say that a Banach space
$X$ satisfies Grothedieck's 
 theorem (in short G.T.) 
if every operator  from $X$
into a Hilbert  space is 1-summing.
\par
We will show that $H^1$-projective
  spaces  are G.T.  spaces of
 cotype 2.
  To this end, we will use the 
following  result [{\bf P1}, Theorem
6.8] .
\par
\proclaim Proposition 4.2. 
A Banach space $X$ is a
$\hbox{G.T.}$  
space of cotype 2, if and only if,
  there are a  metric surjection
$~\sigma : \ell^1(I)\rightarrow X$
 and a constant $c$ such that for
  every 
$x_1,...,x_n $ in $X$  
there are 
$\tilde x_1,...,\tilde x_n $ in
 $ \ell^1(I)$  
such  that  
$$
\forall k\leq n,\quad\sigma (\tilde
 x_k)=x_k\quad\hbox{ and }\quad
\ee \N{\sum_{k=1}^n \eps_k \tilde
x_k }_{\ell^1(I)}  \leq  c 
\ee \N{\sum_{k=1}^n \eps_k x_k }_X
.\eqno (4.1) $$
\par
Let us make some comments. 
A result  from [{\bf P2}] says that
  there exist  numerical  constants
$\alpha, \beta > 0$  
such that  for  every complex  
 Banach space 
 $X$, the following holds:
$\forall x_1,...,x_n  \in X$ 
$$
\beta \ee \N{\sum_{k=1}^n \eps_k x_k 
}_X  \leq 
\int_\tore \N{\sum_{k=1}^n
e^{i3^k\theta} x_k  }_X  \,dm(\theta)
 \leq \alpha  
 \ee \N{\sum_{k=1}^n \eps_k x_k }_X $$
So we can replace (4.1) by 
$$
\N{\sum_{k=1}^n
e^{i3^k\theta} \tilde x_k  }_{\tilde
H^1(\ell^1(I))} \leq c 
\N{\sum_{k=1}^n
e^{i3^k\theta} x_k  }_{\tilde
H^1(X)} \eqno (4.1)'  $$
\par
On the other hand, let $X$ be 
a complex Banach space and 
let $S$ be a subspace  of some
$L^1$.  
Denote  by $ \widetilde S(X) $ the 
  closed  subspace of $ L^1(X)$
generated by $ S \otimes X$.  
Clearly there is a canonical 
embedding   $ S \hat\otimes X 
\hookrightarrow \widetilde S(X)$.
We say that $X$ is $S$-projective 
 if this map is an isomorphism.
Let Rad$^1$ be the closed span of 
$\{ \eps_n : n\geq 1 \} $  
in  $ L^1(\mu)$.
\par
Now, if we put  $ S= H^1$, we obtain
$H^1$-projective  
spaces whereas; if $S~$ is either
Rad$^1$  or  equivalently
$\overline{\hbox{span}}\{e^{i3^n(.)}:
n\geq 1 \} $
then,  by Proposition 4.2  we obtain
precisely G.T. spaces
of cotype 2.
\par
These remarks explain  the
analogy between  $H^1$-projective  
spaces  and G.T. spaces  of cotype 2, and are also  behind  the proof of
the following  theorem .
 \bg
\proclaim Theorem 4.3.
Every $ H^1$-projective  
Banach space  is  a G.T. space  of 
cotype 2.
\par
\noindent Proof :  
Our proof does  not  use the  result
of  [{\bf P2}] but, instead, the  
following   well-known lemma,  
whose proof is elementary and 
 omitted. 
\par
\proclaim Lemma. 
If  $f$ is a continuous  function,
$f : \tore^{p+1} \rightarrow \comp ~$
then
$$
\lim_{n\rightarrow \infty}
\int_\tore f(\theta,n\theta,...,
n^p\theta) \,dm(\theta) =
\int_{\tore^{p+1}} f(\theta_0,...,
\theta_p) \,
dm(\theta_0)...dm(\theta_p) .  $$
\par
 Let $X~$ be an $H^1$-projective  
Banach space. 
Then there exists 
$\sigma :\ell^1(J)\rightarrow x ~$
such that 
$\tilde \sigma :\widetilde H^1(\ell^1(J))
\rightarrow \widetilde H^1(X) $  
is an $\eta (X)$-surjection.
\par
Let $ x_1,...,x_m\in X$, 
using the  lemma we  find $n>0$
such  that 
$$
\int_\tore \N{\sum_{k=1}^m 
e^{in^k\theta} ~x_k}_X \,dm(\theta)
\leq 2 \int_{\tore^m}
\N{\sum_{k=1}^m 
e^{i\theta_k} ~x_k}_X
\,dm(\theta_1)...dm(\theta_m)  \eqno
(4.2) $$
But the function $ f(\theta) =
\sum_1^m e^{in^k\theta}~x_k~$ is
in $\widetilde H^1(X)$  so we can find
$ h = \sum_0^\infty  e^{ip\theta}~
h_p \in \widetilde H^1(\ell^1(J)) $
  with 
$$
\tilde \sigma (h) = f\quad \hbox{ and}
\quad\N{h}_{\widetilde H^1(\ell^1(J)) }
\leq \eta (X) \N{f}_{\widetilde H^1(X)}
\eqno (4.3) $$
For a fixed $j\in J~$,  we have 
by Paley's  inequality 
[{\bf Z}, II.p.121 ]
$$
\left (
\sum_{k=1}^m \abs{ h_{n^k}(j)}^2
\right )^{1\over 2} \leq 2 
\int_\tore \abs{\sum_{p=0}^\infty 
h_p(j) e^{ip\theta }} \, dm(\theta) 
\eqno (4.4) $$
On the other hand, we have 
$$
\ee \abs{ \sum_{k=1}^m \eps_k
  h_{n^k}(j) } \leq \left ( \ee
\abs{ \sum_{k=1}^m \eps_k 
 h_{n^k}(j) }^2 \right )^{1\over 2} 
= \left (
\sum_{k=1}^m \abs{ h_{n^k}(j)}^2
\right )^{1\over 2}
\eqno (4.5) $$
Putting  together (4.4)  and (4.5),  
and taking the  sum  over  all 
$j \in J$; we get 
$$
\ee\N{\sum_{k=1}^m
 \eps_k h_{n^k} }_{\ell^1(J)}
\leq 2 \N{h}_{\widetilde H^1(\ell^1(J)) }
\eqno (4.6) $$
Let us define  $ \tilde x_k = h_{n^k}
\in \ell^1(J) $. 
Using (4.2), (4.3)  and (4.6),   we
obtain 
$\forall\,  k\leq m,\quad\sigma (\tilde
x_k) = x_k$ and 
$$
\ee\N{\sum_{k=1}^m
 \eps_k \tilde x_k }_{\ell^1(J)}  
\leq 4 \eta (X) 
\ee \N{\sum_{k=1}^m
 e^{i\theta_k}  x_k
}_X   \eqno (4.7)  $$
But in general we always have 
$$
\ee \N{\sum_{k=1}^m e^{i\theta_k}  x_k
}_X  \leq 2 
\ee\N{\sum_{k=1}^m
 \eps_k  x_k }_X . $$
So finally, we have proved that 
for   every 
$x_1,...,x_m ~$ in $X$  
there are 
$\tilde x_1,...,\tilde x_m$ in
 $ \ell^1(J)$  
such  that  
$$
\forall k\leq m,\quad \sigma (\tilde
 x_k)=x_k\quad\hbox{ and }\quad
\ee \N{\sum_{k=1}^m \eps_k \tilde
x_k }_{\ell^1(I)}  \leq  8 \eta (X) 
\ee \N{\sum_{k=1}^m \eps_k x_k }_X
 $$ 
This, by Proposition 4.2,
 proves that $X$ is a G.T. space
of cotype 2.\cqfd
\bg
\noindent{\it Remark .} 
One  can  deduce Proposition 3.1
from Theorem 3.5 and Theorem 4.3;  
this yields  a proof  that   does 
not make  use  of the result  of
[{\bf GM}]. 
\bg
Using  deeper  methods,  
we  were  able  in [{\bf K}] to
prove the following theorem.\bg
\proclaim Theorem 4.5. Every
$ H^1$-projective space
$X$ can be isometrically  
embedded  in an $ H^1$-projective   
space $Y$ satisfying \sn
\item{$i.$} $Y \hat\otimes 
Y = Y\check \otimes Y$.
\item{$ii.$} $Y^\ast$ is  a G.T. space of cotype 2.
\item{} where $ Y\check \otimes Y $ is 
the  injective  tensor product of  $Y$ by $Y$.
\bg
\noindent{\bf 5. H${}^1$-projectivity,
 and ultraproducts}
\bg
The main result of this  section
is Proposition 5.4, which asserts  
that the class  of $ H^1$-projective
 Banach spaces  is  closed under the 
formation of ultraproducts .  
 We now give the definition.  Let
$(E_i)_{i \in I}  $
 be a familly  of Banach spaces .  
Consider  the space 
$ \ell^\infty(I,(E_i)_{i \in
I})  $ 
of families 
$(x_i)_{i \in I}  $ 
with $x_i \in E_i~
~(i \in I )$   and  
$$
\N{(x_i)_{i \in I}} =
\hbox { Sup }\{
\N{x_i}_{E_i} : i
\in I \}  < +\infty .$$
$ \ell^\infty(I,(E_i)_{i\in
I})  $ equipped with this norm is
a Banah  space .  
Let $ {\cal U} $ be an ultrafilter
on $ I $ and let $ N_{\cal U} $ 
be the subset  of all those families
$(x_i)_{i \in I} \in
\ell^\infty(I,(E_i)_{i \in
I})  $  with 
$$
\lim_{\cal U}\N{ x_i }_{E_i} = 0 .$$
Obviously, $ N_{\cal U} $  
is a linear subspace of  
$\ell^\infty(I,(E_i)_{i \in
I})  $ 
and  it follows by a standard argument
  that  $ N_{\cal U} $  is closed . 
Now we are ready to recall  the 
following definition.
\par
\noindent{\bf Definition 5.1}.  The  ultraproduct
$(\prod E_i)/{\cal U} $
of the family of Banach  spaces
$(E_i)_{i \in I}  $  
with respect to the ultrafilter 
$ {\cal U} $ is the quotient space
$\ell^\infty(I,(E_i)_{i \in
I}) / N_{\cal U} $ 
equipped with  the canonical quotient
norm . 
\par
For more about ultraproducts,  we
refer the  reader to [{\bf H}].
\par
We fix the following  notation. 
If $X$ is a complex Banach  space,
let $ {\cal P}_n(X) $ denotes the set
of analytic $X$-valued  polynomials  
of degre $\leq n$ :
$$ 
{\cal P}_n(X) =\left \{\sum_{k=0}^n
z^k a_k : z \in \tore \hbox { and }
a_k \in X \right \}    $$
For  a positive  integer r we
denote  by  $ A_r $  
the set $\{ \omega_r^k :
k=0,...,r-1\} $  of $r$th  roots of
1, \ie $ \omega_r = \exp
({{2\pi i}\over r})$.
\proclaim Lemma 5.2. 
Let $X$ be a complex Banach space . 
Then  $\forall n\geq 1, \forall \eps
\in ]0,1[, \forall s \geq \left (
1 + \left [{2\over \eps} \right ]
\right ) n $, and $\forall f \in
{\cal P}_n(X) $  we have
$$
(1-\eps)\int_\tore
\N{f(e^{i\theta})}_X \,dm(\theta) \leq
{1\over s}\sum_{\omega \in A_s}
\N{f(\omega)}_X \leq (1-\eps)^{-1} 
\int_\tore
\N{f(e^{i\theta})}_X \,dm(\theta) $$
\par
The  proof of this lemma can be 
found in [{\bf Z},Ch.X], we include it
for the convenience of the reader .
\par
\noindent Proof :  
Note first that if 
$ g = \sum_{-m}^m z^k a_k  $ 
with  $ a_k \in X $,  then
$$
\forall s>m,\quad\int_\tore
g(e^{it}) \, dm(t) = 
{1\over s}\sum_{\omega \in A_s}
g(\omega ) .\eqno (5.1) $$
Let now  $f$ be an element  of 
${\cal P}_n(X)$, then for every
$r > 1$ 
$$f(e^{i\theta}) = f \ast
V_{n,r}(\theta) = {1\over {2\pi}}
\int_0^{2\pi} f(e^{it})
V_{n,r}(\theta-t) \, dt,\eqno (5.2) $$
and  using (5.1)  we  obtain that
$$
\forall s\geq n(r+1),\quad
f(e^{i\theta}) = {1\over s}
\sum_{k=0}^{s-1} f(e^{{{2\pi i}\over
s}k}) V_{n,r}(\theta - {{2\pi}\over
s}k). \eqno (5.3) $$
By (5.2), for every  $ s \geq n(r+1) $
we have
$$
{1\over s}
\sum_{k=0}^{s-1} \N{ f(e^{{{2\pi
i}\over s}k}) }_X  \leq {1 \over{2\pi}}
\int_0^{2\pi}\N{f(e^{it})}_X\left
( {1\over s}\sum_{k=0}^{s-1}
\abs{ V_{n,r}({{2\pi}\over s}k - t)}
\right ) \, dt \eqno (5.4) $$
but from the definition of the
De La Vall\'ee-Poussin kernel we have
$$
{1\over s}\sum_{k=0}^{s-1}
\abs{ V_{n,r}({{2\pi}\over s}k - t)}
\leq {1\over{r-1}} \left (
{r\over s}\sum_{k=0}^{s-1}
 K_{rn}({{2\pi}\over s}k - t) +
{1\over s}\sum_{k=0}^{s-1}
 K_n({{2\pi}\over s}k - t)\right ) $$
so using (5.1), we obtain easily that
$$
{1\over s}\sum_{k=0}^{s-1}
\abs{ V_{n,r}({{2\pi}\over s}k - t)}
\leq {{r+1}\over{r-1}} $$
and replacing this majorization  
in (5.4), we get that for every 
$ s \geq n(r+1) $
$$
{1\over s}
\sum_{\omega \in A_s} \N{
f(\omega) }_X  \leq {{r+1}\over{r-1}}
{1 \over{2\pi}}
\int_0^{2\pi}\N{f(e^{it})}_X\,dt
\eqno (5.5) $$
\par
On the other hand, starting from
(5.3) we have
$$\eqalign{
\int_\tore\N{f(e^{i\theta})}_X
\,dm(\theta) &\leq {1\over s}
\sum_{k=0}^{s-1} \N{f(e^{{{2\pi
i}\over s}k})}_X 
\int_\tore \abs{V_{n,r}(\theta -
{{2\pi}\over s}k)}\,dm(\theta) \cr
&\leq \N{V_{n,r}}_1 \left ({1\over s}
\sum_{k=0}^{s-1} \N{f(e^{{{2\pi
i}\over s}k})}_X \right ) \cr
&\leq {{r+1}\over{r-1}} {1\over s}
\sum_{\omega \in A_s} \N{
f(\omega) }_X .}\eqno (5.6)$$
It is then sufficient to  choose  
$r = \left [{2\over \eps}\right ]$
  in (5.5) and (5.6).\cqfd
\bg
The following caracterization 
of $ H^1$-projectivity 
 is  more adapted to ultraproducts.
\par
\proclaim Proposition 5.3. 
For $ n\geq 2 $, let $ B_n = A_{12n}$
 and let $X$ be a complex Banach
space. The following assertions are
equivalent.\sn 
\item{$i.$} $X$ is $H^1$-projective.
\item{$ii.$} There exist a constant $c$ and  a metric surjection.
$\sigma : \ell^1(I) \rightarrow X$ such that 
$$
\forall n\geq 1,\quad\forall f \in 
{\cal P}_n(X),\quad\exists  g \in
{\cal P}_{3n}(\ell^1(I)) ~:~
\tilde \sigma (g) = f\quad\hbox{ and }\quad
 \sum_{\omega \in B_n} \N{
g(\omega) }_{\ell^1(I)} \leq c 
\sum_{\omega \in B_n} \N{
f(\omega) }_X . $$

\par
Note that by Lemma 5.2, for every
complex Banach  space $Y$, we have
$$ \forall n\geq 1, \quad\forall  
h \in {\cal P}_{3n}(Y), \quad
{1\over 3} \int_\tore
\N{h(e^{i\theta})}_Y \,dm(\theta) \leq
{1\over {\abs{B_n}}}\sum_{\omega \in
B_n} \N{h(\omega)}_Y \leq
3  \int_\tore
\N{h(e^{i\theta})}_Y \,dm(\theta). $$
\par
Using this remark, the proof of 
Proposition  5.3  becomes  very  easy
and is left as an exercise for the
reader, in particular we can take
$c = 18 \eta (X) $ in $i.
\Rightarrow ii.$ .
\par
\proclaim Theorem 5.4. 
Let $(X_i)_{i \in I}  $ 
be a family of $ H^1$-projective   
Banach spaces such that 
$$ 
\eta = \hbox{ Sup } \{ \eta
(X_i) : i \in I \} <
+\infty, $$
then for every  ultrafilter 
 ${\cal U}$  
on $I$ the ultraproduct 
$(\prod X_i)/{\cal U}$  
is $H^1$-projective.
\par
\noindent Proof : We know that for each $ i
\in I $ there  exists a metric
surjection 
$ \sigma_i : \ell^1(J_i)
\rightarrow X_i$  such that 
$$
\forall n\geq 1,~\forall f \in 
{\cal P}_n(X_i),~\exists  g
\in {\cal P}_{3n}(\ell^1(J_i))~:~
\tilde \sigma_i (g) = f~ \hbox{
and }~
 \sum_{\omega \in B_n} \N{
g(\omega) }_{\ell^1(J_i)} \leq
c_i  \sum_{\omega \in B_n} \N{
f(\omega) }_{X_i} \eqno (5.7)$$ 
Moreover, we have $ c_i \leq K
\eta $,  where $ K$ is some numerical
constant. \par Let $ L=(\prod \ell^1
(J_i))/{\cal U} $, 
we  know (see [{\bf H}]) that $L$
is an  $L^1$-space.
\par
Let $ q : \prod X_i \rightarrow
X =(\prod X_i)/{\cal U} $ 
(resp. $ \bar q : \prod \ell^1(J_i) 
\rightarrow L $) 
be the canonical  quotient  map, and
let $ \sigma =(\prod
\sigma_i)/{\cal U} :L
\rightarrow X $  be defind by $
\sigma(\bar q((x_i)_{i \in I})) =
q((\sigma_i(x_i))_{i \in
I} ) $ for all 
$(x_i)_{i \in I} \in 
\prod X_i $.  
It is immediate  
to see that  $\sigma $  
is  a  metric  surjection.
\par
Fix  now  $ n\geq 1 $  
and  $ f = \sum_0^n z^k a_k \in 
{\cal P}_n(X) $. 
By definition  of the ultraproduct, 
for each $ k \in \{ 0,...,n\} $  
we  can find  $ (a_k^i)_{i
\in I} \in
\ell^\infty(I,(X_i)_{i \in
I}) $  which  represents $ a_k \in X.$
\par
For $ i \in I $,  let $f_i $ 
 denotes  the  element 
$\sum_0^n z^k a_k^i \in 
{\cal P}_n(X_i)   $  
then for  each  $z \in \tore $  
we have $ f(z) = q((f_i(z))_
{i \in I} ) $.  
Consequently, for  every $\omega \in
B_n $  there exists $ I_\omega  \in
{\cal U} $ such that  
$$ \forall  i  \in I_\omega,\quad
\N{f_i (\omega) }_{X_i }
\leq 2 
\N{f(\omega) }_X, $$
since $\lim_{\cal U} 
\N{f_i (\omega) }_{X_i } = 
\N{f(\omega) }_X$. Define then 
$ I_n = \cap_{\omega \in B_n }
I_\omega \in {\cal U} $,  we obtain
$$
\forall  i  \in I_n ,\quad\forall 
\omega \in B_n ,\quad
\N{f_i (\omega) }_{X_i }
\leq 2 
\N{f(\omega) }_X \eqno (5.8) $$
\par 
Let us now define $(b_k^i )_{ i  \in
I} $  by $ b_k^i  = a_k^ i  $
 if $  i  \in I_n $, and   
$ b_k^i  = 0 $
 if $  i  \not\in I_n $. Put
$ g_i  = \sum_0^n z^k b_k^ i 
\in  {\cal P}_n(X_i )   $.  
Clearly  for each $ k \leq n $  
we have $ a_k =
q((b_k^i )_{ i  \in I} ) $, 
and by (5.8) we get 
$$
\forall  i  \in I ,\quad\forall 
\omega \in B_n,\quad
\N{g_ i (\omega) }_{X_i }
\leq 2 
\N{f(\omega) }_X \eqno (5.9) $$
\par
Now using (5.7) we can find for 
each $ i \in I $ an element
$ h_i  = \sum_0^{3n} z^k t_k^i 
\in {\cal P}_{3n}(\ell^1(J_i ))
$  
such  that  $ \tilde
\sigma_i (h_i ) = g_i  $  
and
$$
\sum_{\omega \in B_n} \N{
h_ i (\omega) }_{\ell^1(J_i )}
\leq K \eta  \sum_{\omega \in B_n}
\N{ g_i (\omega) }_{X_i } 
\leq 2 K \eta \sum_{\omega \in B_n}
\N{ f(\omega) }_X, \eqno (5.10) $$
where the last inequality  comes
 from (5.9).\par
On the other hand,  it is clear  
that $ (t_k^i )_{ i  \in I} 
\in  \ell^\infty (I,(
\ell^1(J_i ))_{ i  \in I} ) $  
therefore  if we put $
t_k =\bar q((t_k^i )_{ i  \in I})
 \in L $  and  consider  
$h = \sum_0^{3n} z^k t_k
\in {\cal P}_{3n}(L)  $,  
we obtain $\tilde\sigma (h) = f $  
and from (5.10), we deduce that
$$
\sum_{\omega \in B_n} \N{
h(\omega) }_L
=\lim_{\cal U}  \sum_{\omega \in B_n}
\N{ h_ i (\omega) }_{\ell
^1(J_ i )} 
\leq 2 K \eta \sum_{\omega \in B_n}
\N{ f(\omega) }_X  . $$
the result follows  now by another
use of Proposition 5.3. \cqfd
\bg
\noindent{\it
Remark. } 
In fact we can show  that 
$$
\eta \left ((\prod X_ i )/{\cal U}
\right ) \leq \eta = \hbox{ Sup }
 \{ \eta (X_i ) :  i  \in I \}
$$
by  refining  upon  the  preceding
arguments.
\par
Using this remark and the fact that
$ X^{\ast\ast}$ is 1-complemented  
in some ultrapower $ (\prod X)/{\cal
U} $  of $X$, we  obtain another
proof of Corollary 2.6.
\par
\bg
\noindent\bf{6. Exemples  and concluding   
remarks }
\par 
\bg
\rm
{\bf a}. 
Every ${\cal L}^1$-space  is 
$ H^1$-projective.  
This  can be seen directly, or by
using Corollary (2.6)  and the
following fact from [{\bf LR}] :  
the bidual of  an  ${\cal L}^1$-space
 is isomorphic  to  a complemented
subspace  of an $L^1$-space.
\par
 This statement  does  not  hold  for
the non-commutative   
analogues of  $L^1$-spaces    
Indeed, the trace class  operator
   ideal $ C_1 = \ell^2 \hat\otimes
 \ell^2 $ contains a complemented
copy of $ \ell^2 $ therefore  $C_1 $
  is not  $ H^1$-projective. Note
that $ C_1 $, fails  the U.H.M.D.  
property  as  it is shown  in 
 [{\bf HP}]. 
\par
{\bf b}.  
It is known (\cf [{\bf B}] )  that
 $L^1/H^1$ is a G.T.   space  of 
cotype 2,  but  by a well-known  
conterexample  we know also  that
 $L^1/H^1$  fails the analytic
  Radon-Nikodym  property. So the
convese of Theorem 4.3  does  not 
hold.
\par
{\bf c}.  
 The  next  theorem  provides us  with
   non-trivial   examples  of 
$H^1$-projective    
spaces, its proof  was  implicit in
[{\bf BD}],   we sketch  it  for  
the convenience of   the  
reader. 
\par
\proclaim Theorem 6.1. 
Let $ Y $ be  a  reflexive  subspace  
of  $ L^1(\Omega,{\cal A},\mu) $,  
then $ L^1(\Omega,{\cal A},\mu) / Y$
 is  $ H^1$-projective.
\par
\noindent Proof : 
The proof is based on the following 
two facts :
\par
I. If $(\Omega,\mu )$  
is  a measure space,  then for every
$ \alpha \in ] 0, 1[ $ 
 the  Riesz  projection $\Re $  
is  a continuous  operator  
from $L^2\left (\tore,dm;
L^1(\Omega,\mu) \right ) $
into
 $L^2\left (\tore,dm;
L^\alpha(\Omega,\mu) \right ) $ .
 See [{\bf BD}]. 
\par
II. Let $Y$ be a  reflexive  
subspace of $ L^1(\Omega ,\mu) $, 
then,  modulo  a change of density  
if necessary,  the normes $\N{.}_1$ 
and $\N{.}_p$  are   equivalent  
on  $Y$  for some 
$ p > 1 $. 
 See [{\bf R}]. 
\par
Using Fact II  and 
 H\"older's
inequality  
we see  easily  that   
$$
\forall\,y \in Y,\qquad
\N{y}_{L^1(\mu)} \leq c_1
\N{y}_{L^{1/2}(\mu)}  \eqno (6.1) $$
Let $ q: L^1(\Omega ,\mu) 
\rightarrow L^1(\Omega ,\mu)/Y $ 
be  the  quotient  map, and 
consider an element $ f \in \widetilde H^2(L^1(\Omega,\mu)/Y)  $  
such  that  $ \N{ f }_{\tilde
H^2(L^1/Y)} < 1 $. 
We can find $ g \in L^2\left (\tore,dm;
L^1(\Omega ,\mu) \right ) $  
satisfying 
$$ 
\tilde q(g) = f \quad\hbox{ and }\quad
\N{g}_{L^2(L^1)} < 1.\eqno (6.2) $$
By Fact I, the negative Riesz 
projection $ \Re_- $ is 
bounded   from $L^2\left (\tore,dm;
L^1(\Omega ,\mu) \right ) $
into
 $L^2\left (\tore,dm;
L^{1/2}(\Omega ,\mu) \right)$. 
  So  using (6.1), we get 
$$
\N{\Re_-(g)}_{L^2(L^{1/2})}  \leq c_2. \eqno (6.3)   $$
Since $ \tilde q(\Re_-(g)) = \Re_-
(\tilde q(g)) = 0 $,  we see that  
$ \Re_-(g) $  takes  its  values 
in  $Y$,   so by  (6.1) and  (6.3)
$$ 
\N{\Re_-(g)}_{L^2(L^1)}  \leq c_1c_2.$$
Therefore, if we put $ h = g - \Re_-(g)
\in \widetilde H^2(L^1(\Omega,\mu))$
we have
$$\N{h}_{\tilde
H^2(L^1(\mu)) } \leq (1+ c_1c_2 )=c_3.
$$
Consequently, we  have  proved 
that   for  every  $ f \in \tilde
H^2(L^1(\Omega,\mu)/Y)  $  there
exists  $ h  
\in \widetilde H^2(L^1(\Omega,\mu))$  
such  that $ \tilde q(h) =f $   and
$$\N{h}_{\tilde
H^2(L^1(\mu)) } \leq c_3
 \N{ f }_{\tilde
H^2(L^1/Y)}. 
$$
the result  follows now from
Proposition 1.1 and Proposition 2.2.\cqfd
\par
{\bf d}. 
It might  be of  some interest  to
have  estimations  for the
$ H^1$-projectivity constant of
finite-dimensional spaces.
\par
Recall that the Banach-Mazur distance
$~d(E,F)$ between  two
 finite-dimensional  spaces $~E,~F $,
with $\dim E=\dim F $ is equal  to 
$ \inf \{\N{T}\N{T^{-1}}\} $, 
where  the infimum  is taken over  
all  isomorphisms  $ T: E\rightarrow 
 F $.
\par
It  is immediate  
 that $ \eta (E) \leq d(E,F) \eta
(F) $  
for  every  pair  of
finite-dimensional spaces $ E, F$  
 of the same dimension. 
Therefore, we always have  
$$
\eta (E) \leq d(E,\ell_n^1 ) .\eqno (6.4)  $$
In fact we can  also prove the 
following.
\par
\proclaim Proposition 6.2. 
For every $ p \in [1,+\infty] $ 
and  $ n\geq 1 $ we have the
following 
$$
{1\over 2}. n^{1-1/r} \leq 
 \eta (\ell_n^p) \leq
d(\ell_n^p,\ell_n^1) \leq  n^{1-1/r} 
\eqno (6.5)$$
where  $ r = \min (2,p) $.
\par
\noindent Proof : 
The inequality $
d(\ell_n^p,\ell_n^1) \leq  n^{1-1/r} $
  is  classical. To see it
elementarily  take for $ T: \ell_n^1
\rightarrow \ell_n^p $  the operator
 defined  by  $ T(x) = x $ if 
$ p \in [1,2] $  
and  by $ T(x) = W_n(x) $ 
(where $W_n $  is the  matrix
$(a_{jk})_{1\leq j,k \leq n} $  
with $ a_{jk} = \exp [{{2\pi i}\over
n } (j-1)(k-1)] $) 
if $ p \in [2,+\infty] $.
\par
Using  Paley's inequality, we  
know that  the operator 
$ u : H^1\rightarrow \ell_n^2 $ 
 defined  by  $
 u(h) = (\hat h(3),...,\hat h(3^n)) $
  is of norm smaller than 2 . 
Therefore, 
$$
\forall q \in [1,+\infty] ,\quad
\N{u}_{{\cal L}(H^1,\ell_n^q)}
 \leq 2 \max (1, n^{1/q - 1/2} ) .
\eqno (6.6) $$
\par
Let $ e_1,...,e_n $  
be the canonical basis  of $ \ell_n^p
$
 and consider $ f \in  \widetilde H^1
(\ell_n^p))$ defined by 
$ f(\theta) = \sum_{k=1}^n
e^{i3^k\theta}~e_k $,  we have
$$
\N{f}_{\widetilde H^1(\ell_n^p)} =
n^{1/p} . \eqno (6.7) $$
 Assume that $ {1\over p} + {1\over q}
= 1 $. Using the duality 
$ (H^1\hat\otimes \ell_n^p)^\ast =
 {\cal L} ( H^1,\ell_n^q) $  
 and the  fact that $$
\langle f, u\rangle = \sum_{k=1}^n
\langle e_k, u ( e^{i3^k(.)}) \rangle 
 = n $$ we obtain 
$$
n \leq \N{f}_{H^1\hat\otimes \ell_n^p}
~\N{u}_{ {\cal L} ( H^1,\ell_n^q)} 
\leq 2 \max (1, n^{1/q - 1/2} )
\N{f}_{H^1\hat\otimes \ell_n^p}, $$
the last inequality comes from (6.6). 
On the other hand, if  we  use (6.7) 
 and  the definition  of
$ H^1$-projectivity,  
we deduce easily that  
$$
n \leq 2 \max (1, n^{1/q - 1/2}) 
\eta (\ell_n^p)  n^{1/p}   $$
which yields the first inequality 
in  (6.5) and proves the result.\cqfd
\par
{\bf e}.
By the results of section 4. and a
theorem of Pisier [{\bf P1},Ch 4],  
we see  that an $ H^1$-projective  
space $X$ having  the 
approximation 
 property  and  such  that 
 $X^\ast$  is also  $
H^1$-projective must  be 
finite-dimensional . So we will end 
our  discussion  by asking  the
following  question : Can one remove the approximation  
hypothesis   
from the preceding statement ? 
\par 
\bg
{\sc Acknowledgements :}
The auther  thanks Professor Gilles
 Pisier for pointing  out this subject
to him  and for many stimulating 
discussions  during  the preparation  
of  this  paper. 
\par
\bg
\centerline{\ninerm REFERENCES }
\par
\ref B & J.Bourgain & New Banach 
space properties of the disc algebra 
and $H^\infty $, &Acta Math. & 
152 (1984), 1-48.

\ref BD & J.  Bourgain and W.J. Davis
&  Martingales transforms and complex
uniform convexity,& Trans. Amer.
Math. Soc.& 294 (1986),501-515.

\ref CW & R. Coifman and G. Weiss
& Extensions of Hardy spaces and their
  use in analysis, & Bul. Amer. Math.
Soc.& 83 (1977).

\ref D & D.W. Dean & The equation,
${\cal L}(E,X{\ast\ast}) = 
{\cal L}(E,X)^{\ast\ast} $ and the 
principle of local reflexivity, &
Proc. Amer. Math. Soc. & 40 (1973), 146-148.

\ref E1 & G. Edgar && Complex
martingale convergence, & Lecture
Notes in Mathematics, 1166,
Springer-Verlag (1985), 38-59.

\ref E2 & G. Edgar & Analytic
martingale convergence, & J. Funct.
Anal.& 69 (1986), 268-280.

\ref Ga & D.J.H. Garling & On
martingales with values in complex
Banach spaces, & Proc. Cambridge
Philos. Soc. & 104 (1988), 399-406.

\ref GLM & N. Ghoussoub, J.
Lindenstrauss and B. Maurey &
Analytic martingales and
plurisubharmonic  barriers in
complex Banach spaces, &
Contemporary Mathematics, &Vol.85, Amer. Math. Soc., Providence, RI, (1989). 111-130

\ref GM & N. Ghoussoub and B. Maurey
& Plurisubharmonic martingales  and 
barriers in complex quasi-Banach
spaces,&Ann. Inst. Fourrier& 39 (1989), 1007-1060.

\ref G1 & A. Grothendieck &
R\'esum\'e de la th\'eorie  m\'etric
des produits tensoriels topologiques,&
Bol. Soc. Math. S\~ao Paulo.& 8 (1956)
 1-79.

\ref G2 & A. Grothendieck &
Produits tensoriels topologiques,&
Memoires Amer. Math. Soc.& 16 (1955).

\ref HP & U. Haagerup and G. Pisier &
Factorization of analytic  functions 
with values in non-commutative
$L^1$-spaces, & Canad. J. Math. & 41 (1989), 882-906.

\ref H & S. Heinrich & Ultraproducts 
in Banach space theory, & J. Reine
Angew. Math.& 313 (1980), 72-104.

\ref K & O. Kouba & L'application 
canonique $J: \widetilde H^2(X)
\hat\otimes \widetilde H^2(X) \rightarrow
\widetilde H^1(X\hat\otimes X) $ n'est 
pas surjective en g\'en\'eral.&
C.R. Acad. Sci. Paris t.307, S\'erie 
I, & (1988), 949-953.

\ref LR & J. Lindenstrauss and H.P.
Rosenthal & The ${\cal L}_p $ spaces, &
Is. J. Math. & 7 (1969), 325-34.

\ref LT & J. Lindenstrauss and L.
Tzafriri && Classical Banach spaces
I,II ,& Springer-Verlag (1977).

\ref P1 & G. Pisier && Factorization 
of  linear operators and geometry of
Banach spaces, & CBMS No 60 , A.M.S.
Providence (1987) .

\ref P2 & G. Pisier & Les
in\'egalit\'es de Khintchine-Kahane 
d'apr\`es C. Borell. S\'eminaire
sur la g\'eometrie des espaces de
Banach (1977-1978), & Exp. No 7,
Ecole Polytec. Palaiseau & (1978).

\ref P3 & G. Pisier & Factoriztion
of  operator valued analytic
functions, & Advances in Math.& 93 No 1, (1992) 61-125.
 
\ref R & H.P. Rosenthal & On
subspaces of $ L_p $, & Ann. of Math.&
(2)97 (1973), 344-373.

\ref X & Q. Xu & In\'egalit\'es pour 
les martingales de Hardy  et
renormage des espaces quasi-norm\'es, & C.R. Acad. Sci. Paris t. 307
S\'erie I,& (1988), 601-604.

\ref Z & A. Zygmund && Trigonometric
series I, II, & Cambridge University
press the recent edition (1988).\par

\end